\documentclass[11pt ]{article}
\usepackage{amsmath}
\usepackage{graphicx}
\usepackage{hyperref}
\usepackage{lmodern}
\usepackage{amssymb}
\usepackage{booktabs}

\newcommand{\be}{\begin{equation}}
\newcommand{\ee}{\end{equation}}
\newcommand{\bea}{\begin{eqnarray}}
\newcommand{\eea}{\end{eqnarray}}

\renewcommand{\theequation}
{\arabic{section}.\arabic{equation}}

\makeatletter
\def\eqnarray{ \stepcounter{equation} \let\@currentlabel=\theequation
 \global\@eqnswtrue
 \global\@eqcnt\z@
 \tabskip\@centering
 \let\\=\@eqncr
 $$\halign to \displaywidth\bgroup\@eqnsel\hskip\@centering
 $\displaystyle\tabskip\z@{##}$&\global\@eqcnt\@ne
 \hfil$\displaystyle{{}##{}}$\hfil
 &\global\@eqcnt\tw@$\displaystyle\tabskip\z@{##}$\hfil
 \tabskip\@centering&\llap{##}\tabskip\z@\cr}
\makeatother

\makeatletter
\def\@arrayacol{\edef\@preamble{\@preamble \hskip .5\arraycolsep}}
\def\array{\let\@acol\@arrayacol \let\@classz\@arrayclassz
\let\@classiv\@arrayclassiv \let\\\@arraycr\def\@halignto{}\@tabarray}
\makeatother

\makeatletter
\newcounter{subeqncnt}
\def\thesubeqncnt{\alph{subeqncnt}}
\def\subequations{\begingroup%
   \stepcounter{equation}\edef\@tempa{\theequation}%
   \let\c@equation\c@subeqncnt\c@subeqncnt\z@
   \edef\theequation{\@tempa\noexpand\thesubeqncnt}}

\makeatother

\newcommand{\nn}{\nonumber}

\def\CZ {{\cal Z}}

\def\Det{\mathop {\rm Det}}

\begin{document}

\setlength{\baselineskip}{6mm}
\begin{titlepage}
\begin{flushright}

{\tt NRCPS-HE-45-2015} \\

\end{flushright}

\begin{center}
{\Large ~\\{\it

Spectrum and Entropy of C-systems \\
\vspace{1cm}

MIXMAX  Random Number Generator\\
}

}

\vspace{1cm}

 {\sl Konstantin Savvidy$^{@}$   and George Savvidy$^{*}$

 \bigskip
 \centerline{${}$ \sl $^{*}$Institute of Nuclear and Particle Physics}
\centerline{${}$ \sl Demokritos National Research Center, Ag. Paraskevi,  Athens, Greece}
\bigskip

 \centerline{$^{}$ \sl $^{@}$
 College of Science, Nanjing University of Aeronautics and Astronautics}
\centerline{${}$ \sl  Nanjing 211106, China.}

}
\end{center}
\vspace{30pt}

\centerline{{\bf Abstract}}
The  uniformly hyperbolic Anosov  C-systems defined on a torus have very strong instability of
their trajectories, as strong as it can be in principle.
These systems have exponential instability of
all their trajectories and as such have mixing of all orders,
nonzero Kolmogorov entropy and a countable
set of everywhere dense periodic trajectories.
In this paper we are studying the properties of their spectrum and of the entropy.
For a two-parameter family of C-system operators $A(N,s)$,
parametrised by the integers $N$ and $s$, we found the universal limiting form of the spectrum,
the dependence of entropy on  $N$  and the period of its trajectories on a rational
sublattice. One can deduce from this result that the entropy and the periods are sharply
increasing with $N$.
We present a new three-parameter family of C-operators $A(N,s,m)$ and analyse
the dependence of its spectrum and  of the entropy on the parameter $m$.
We are developing our earlier suggestion to use these tuneable  Anosov C-systems
for multipurpose Monte-Carlo simulations. The  MIXMAX family of random number generators
based on Anosov C-systems provide high quality statistical properties, thanks to their large entropy,
have the best combination of speed, reasonable size of the state, tuneable parameters
and availability for implementing the parallelisation.
\vspace{12pt}

\noindent

\end{titlepage}

\section{\it Introduction}
The  uniformly hyperbolic Anosov  C-systems defined on a torus have very strong instability of
their trajectories, as strong as it can be in principle \cite{anosov}.
These systems have exponential instability of
all their trajectories and as such have mixing of all orders,
nonzero Kolmogorov entropy and a countable
set of everywhere dense periodic trajectories\footnote{
D.V. Anosov gave the definition of C-systems in his 
outstanding work \cite{anosov}.
 In order to provide this definition one should use 
such mathematical concepts as the tangent vector bundle, derivative mapping, 
contracting and expanding linear spaces, foliations  and others.  The definition of the 
C-systems \cite{anosov},  of the Kolmogorov entropy \cite{kolmo,kolmo1,sinai3}, the description 
of the properties of 
its periodic trajectories and review of the 
 dynamics can be found in recent article \cite{Savvidy:2015ida}. }.
In this paper we are studying the properties of their spectrum and of their entropy and are
developing our earlier suggestion to use the Anosov C-systems
for Monte-Carlo simulations \cite{yer1986a}.

The particular system chosen for investigation is the one realising linear automorphisms of the unit hypercube in Euclidean space $\mathbb{R}^N$ with coordinates $ (u_1,...,u_N)$\cite{anosov,yer1986a,konstantin,Savvidy:2015ida}:
\be
\label{eq:rec}
u_i^{(k+1)} = \sum_{j=1}^N A_{ij} \, u_j^{(k)} ~~~~~\textrm{mod}~ 1,~~~~~~~~~k=0,1,2,...
\ee
where the components of the vector $u^{(k)}$ are defined as 
$$
u^{(k)}= (u^{(k)}_1,...,u^{(k)}_N).
$$ 
The dynamical system defined here by the integer matrix $A$
should have a determinant equal to one $\Det A =1$.
In order for the automorphisms (\ref{eq:rec})  {\it to fulfill  the Anosov hyperbolicity
C-condition it is necessary
and sufficient that the matrix $A$ has no eigenvalues on the unit circle} \cite{anosov}.
Therefore  the  spectrum $\{ \Lambda = {\lambda_1},...,
\lambda_N \}$ of the matrix $A$ should fulfill the following
two conditions:
\bea\label{mmatrix}
1)~\Det  A=  {\lambda_1}\,{\lambda_2}...{\lambda_N}=1,~~~~~
2)~~\vert {\lambda_i} \vert \neq 1, ~~~\forall ~~i .~~~~~~
\eea
Because the determinant of the matrix $A$ is equal to one,
the Liouville's measure $d\mu = du_1...du_N$ is invariant under the action of $A$.
The inverse matrix $A^{-1}$ is also an integer matrix because $\Det  A=1$.
Therefore $A$ is an automorphism of  the unit hypercube  onto itself.
The conditions (\ref{mmatrix}) on the eigenvalues of the matrix $A$ are  sufficient
to prove that the system represents  an Anosov C-system \cite{anosov} and therefore as such it
also represents a  Kolmogorov K-system \cite{kolmo,kolmo1,sinai3,rokhlin,rokhlin2}
with mixing  of  all orders and of nonzero entropy.

\begin{figure}
 \centering
\includegraphics[width=4cm]{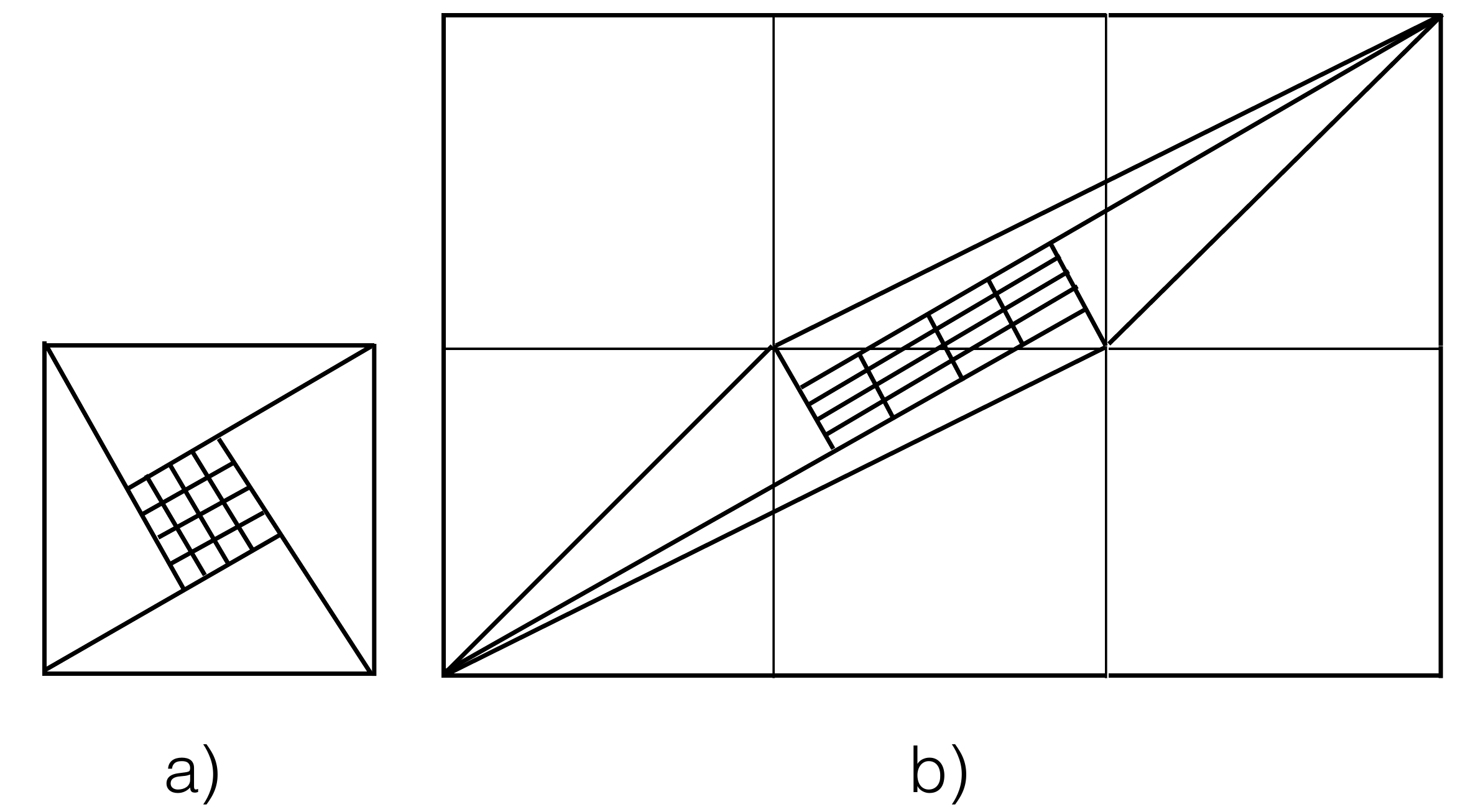}
\centering
\caption{The eigenvectors of the matrix $A$
 $\{ e_{\alpha}  \}$ and  $\{ e_{\beta}  \}$ define
two families of  parallel planes  $\{X_{\alpha} \}$ and $\{Y_{\beta} \}$
which are invariant under the automorphisms  $A$.
The automorphism $A$ is contracting  the distances between points on the
planes belonging to the set $\{X_{\alpha} \}$ and expanding the
distances between points on the planes belonging to the set $\{Y_{\beta} \}$.
The a) depicts the parallel planes of the sets  $\{X_{\alpha} \}$ and $\{Y_{\beta} \}$
and b) depicts their positions after the action of the automorphism  $A$.}
\label{fig1}
\end{figure}

The eigenvalues of the matrix $A$ can always be sorted by increasing absolute value and divided
into the two sets   $\{ \lambda_{\alpha}  \} $ and $\{  \lambda_{\beta }  \} $
with modulus smaller and larger than one:
\bea\label{eigenvalues}
0 <  \vert \lambda_{\alpha} \vert   < 1   & \textrm{ for } \alpha=1...d\nn\\
1 <  \vert \lambda_{\beta} \vert  < \infty & \textrm{ for } \beta=d{+}1...N . 
\eea
There exist two  hyperplanes $ X= \{X_{\alpha} \}$ and $ Y= \{Y_{\beta} \}$
which are  spanned by the
corresponding eigenvectors  $\{ e_{\alpha}  \}$ and  $\{ e_{\beta}  \}$ .
These invariant planes of the matrix $A$,  for which the eigenvalues are outside  and  inside
of the unit circle respectively,  define the expanding  and  contracting invariant spaces, so that
the phase trajectories of the dynamical system (\ref{eq:rec}) are {\it expanding and contracting
under the transformation $A$ at an exponential rate} (see Fig.\ref{fig1}).  The same is true for the inverse evolution which is defined by the matrix $A^{-1}$. For the inverse evolution the
contracting and expanding invariant spaces alternate their role.

The exceptional property of the C-system (\ref{eq:rec}) is that it  has
nonzero entropy and that it is possible to  calculate its value  $h(A)$ in
terms of the eigenvalues of the operator $A$.
The most convenient way to find out the entropy is  to integrate over the whole
phase space the logarithm of the volume expansion rate $\lambda(u)$
of an infinitesimal cube which is embedded  into the expanding foliation
 \cite{anosov,sinai3,rokhlin2,sinai4,gines,Savvidy:2015ida}.
For the automorphisms
(\ref{eq:rec}) the coefficient $\lambda(u)$  does not depend on the phase
space coordinates $u$ and is equal to the product of eigenvalues
$\{  \lambda_{\beta }  \} $  with modulus  larger than one
$
\lambda(u) = \prod_{\beta} \lambda_{\beta}
$, as it was defined in (\ref{eigenvalues}).
Thus for the  Anosov automorphisms   (\ref{eq:rec}) one can calculate the entropy,
which is  equal to the sum:
\be\label{entropyofA}
h(A) = \sum_{\beta   } \ln \vert \lambda_{\beta} \vert.
\ee
{\it This allows to characterise and compare the chaotic properties of different
dynamical systems quantitatively by computing and comparing their entropies. }
As it is obvious from the above formula for the entropy its value depends on the
spectral properties of the evolution operator $A$.  
Also the variety  and richness of the periodic
trajectories of the C-systems  essentially depend on the entropy
\cite{anosov,Savvidy:2015ida,bowen0,bowen}.  Indeed, the number of periodic 
trajectories  $\pi(q)$
of a period $q$ has the form
\be\label{density}
\pi(q)  \sim  e^{q \ h(A)} / q
\ee
and tells that a system with larger entropy $h(A)$ is more densely populated by the
periodic trajectories of the period $q$.
Our aim is to study these characteristics of the C-systems and develop our earlier
suggestion to use the Anosov C-systems for Monte-Carlo simulations
\cite{konstantin,Savvidy:2015ida, yer1986b,mixmaxGalois}.

The earlier publications concerning the application of the modern results of the ergodic theory to
concrete physical systems can be found in \cite{yer1986a,yangmillsmech,Savvidy:1982jk,body}. These articles contain review material as well. 
The recent review articles on random numbers generators 
for the Monte-Carlo simulations can be found in  \cite{pierre1,pierre2}.

In this article we shall explore a two- and  three-parameter family
of matrix C-operators $A(N,s)$ and $A(N,s,m)$ by calculating their
spectrum and the corresponding entropies. This provides us with the knowledge of the
behaviour of the entropy as a function of these parameters and allows to
{\it classify these dynamical systems by the increase of their chaotic-stochastic properties}.
The presence of trajectories of a large period is also associated with the value
of the entropy of a dynamical system  and we shall study the periods
of the above C-system trajectories.
In the second section we shall study all these properties for the two-parameter family
of operators $A(N,s)$. In the third  section the tree-parameter family of
operators $A(N,s,m)$ will be investigated.
In the fourth section we shall consider the application of these
systems to generating random numbers for the Monte-Carlo simulations and other
multidisciplinary purposes.

\section{\it Two-parameters Family of C-operators  $A(N,s)$ }

We shall start with the two-parameter family of operators
introduced in \cite{konstantin}, which are
parametrised by the integers $N$ and $s$. We are interested in investigating
the general properties of the spectrum and the corresponding value of the
Kolmogorov entropy.
The matrix is of the size $N\times N$, its  entries are all integers
 $A_{ij} \in \mathbb{Z}$, and it has the following form \cite{konstantin}:
\be
\label{eq:matrix}
A(N,s) =
   \begin{pmatrix}
      1 & 1 & 1 & 1 & ... &1& 1 \\
      1 & 2 & 1 & 1 & ... &1& 1 \\
      1 & 3{+}s & 2 & 1 & ... &1& 1 \\
      1 & 4 & 3 & 2 &   ... &1& 1 \\
      &&&...&&&\\
      1 & N & N{-}1 &  ~N{-}2 & ... & 3 & 2
   \end{pmatrix}
\ee
The operator is constructed so that its entries are increasing together with the size $N$ of the
operator, and we have a family of operators which are parametrised by the
integers $N$ and $s$.
{\it For any integer values of the parameters $N$ and $s$ of the operator $A(N,s)$, if the condition} (\ref{mmatrix})
{\it is fulfilled, then the operator $A(N,s)$ represents a C-system } \cite{yer1986a}.
In reference \cite{konstantin} and in our Tables \ref{tbl:largeS} and  \ref{tbl:largeN},
some values of $N$ and $s$
were found for  which the C-condition (\ref{mmatrix}) is fulfilled.
As we mentioned above, the chaotic properties of the C-system
operators are quantified  through their spectral distributions,
the value of the Kolmogorov
entropy (\ref{entropyofA}) and periodic trajectories
on a rational sublattice \cite{yer1986a,konstantin,Savvidy:2015ida}.
The data in the Tables \ref{tbl:largeS}, \ref{tbl:largeN}  demonstrates that the entropy
of the operator  $A(N,s)$ strongly  depends on $N$ and $s$
and therefore their  chaotic-stochastic  properties \cite{konstantin} .
\begin{figure}
\begin{center}
\includegraphics[width=5cm]{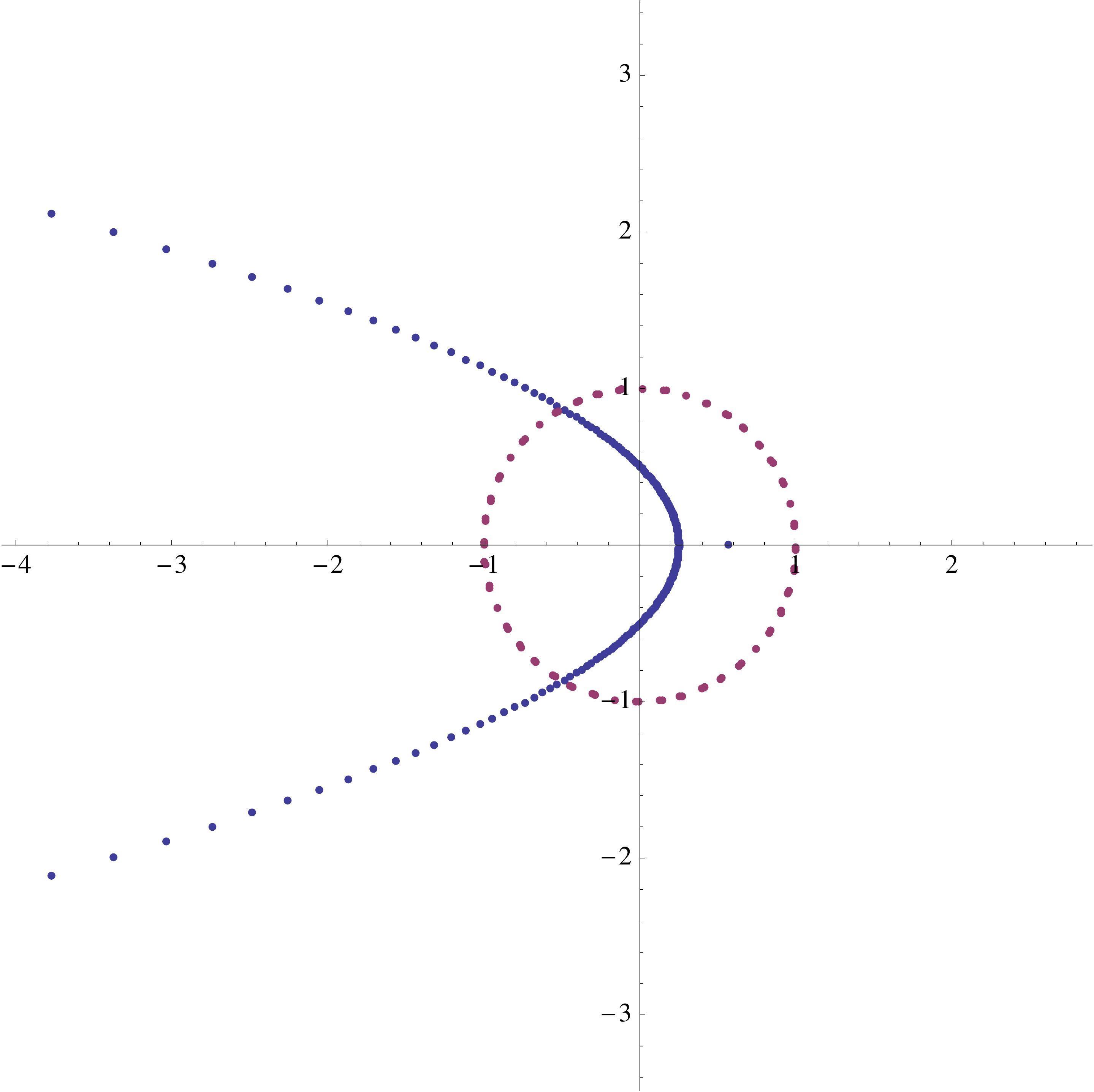}~~~~~~~~~~
\includegraphics[width=5cm]{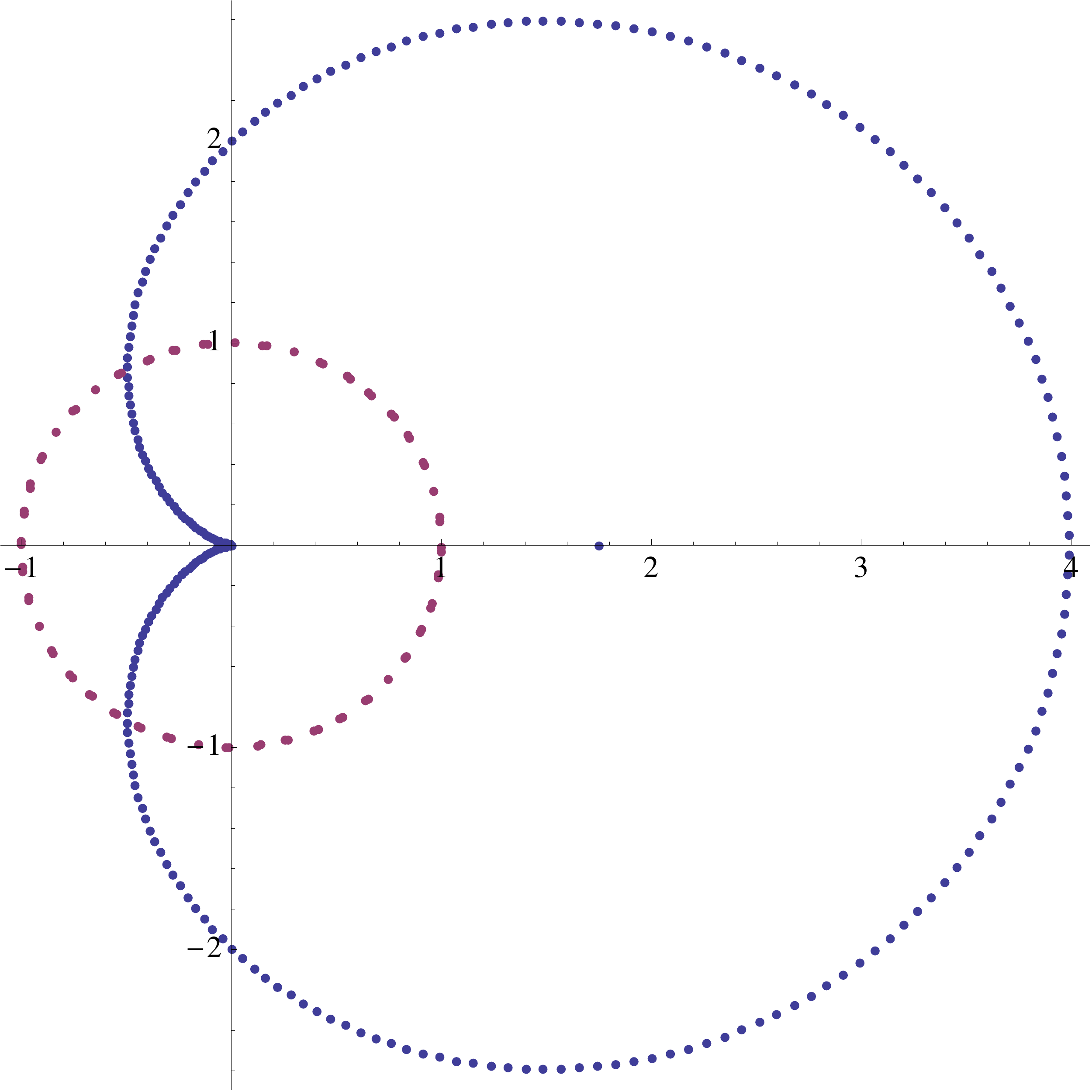}
\caption{
On the left is the distribution of the eigenvalues of the operator $A(N,s)$
and on the right of its  inverse operator $A^{-1}(N,s)$ for the $N=256$ and $s=-1$.
The unit circle is depicted to separate  the eigenvalues inside and outside
the circle in  accordance
with  the formula (\ref{eigenvalues}). The spectrum  
of the operators $A(N,s)$ and $A^{-1}(N,s)$
has two real eigenvalues for even $N$, all the rest
of the eigenvalues are complex. The minimal $\lambda_{min}\approx 0.25$ and the 
maximal $ \lambda_{max}  \approx 3002 $
eigenvalues  are real. The graph on the left  shows  only the small eigenvalues 
which are well approximated by formula (\ref{eq:evs}), the large eigenvalues 
are not depicted.  The spectrum of the inverse operator is 
compact and all eigenvalues can be seen on the right figure.
As $N$ tends to infinity
$N \rightarrow \infty$, the scale and the shape of the complex curves
remain intact, instead,  the eigenvalues are filling these curves more and more densely.
The value $1 /  \lambda_{max}   $ tends to zero and so lies at the stem of
the leaf in the figure on the right.
}
\label{fig2}
\end{center}
\end{figure}

\begin{table}[htbp]
   \centering
   \begin{tabular}{@{} rcccccl @{}} 
      \toprule
      Size & Magic & Entropy & Log of the period $q$ \\
      N    & $s$ &   &  $\approx \log_{10} (q)$  \\ 
      \midrule
        256   & $-1$   & $164.5$           & ${4682}$ \\
        256   & 487013230256099064 & 193.6 & 4682 \\
      \bottomrule
   \end{tabular}
   \caption{Properties of operators $A(N,s$) for large special $s$.  
   The first line is given for comparison, in order to illustrate the improvement of the entropy of the generator for the large $s$ (in the second line).  The period $q$ is defined in (\ref{period1}), (\ref{period}).  }
   \label{tbl:largeS}
\end{table}

\begin{table}[htbp]
   \centering
   \begin{tabular}{@{} lcccrcl @{}} 
      \toprule
      Size & Magic & Entropy & Log of the period q \\
      N    & $s$ &   &  $\approx \log_{10} (q)$  \\ 
      \midrule
        7307   & 0 & 4676.5 & 134158 \\
	20693 & 0 & 13243.5 & 379963 \\
        25087 & 0 & 16055.7 & 460649 \\
 	28883 & 1 & 18485.1 & 530355 \\
	40045 & -3 & 25628.8 & 735321 \\
	44851 & -3 & 28704.6 & 823572 \\
      \bottomrule
   \end{tabular}
   \caption{Table of properties of the operator $A(N,s)$ for large matrix size $N$. The third column is the value of the Kolmogorov entropy, which needs to be greater than about $h(A) \approx 50$ for the generator to be empirically acceptable. Therefore it should not be surprising that for all of these generators the sequence  passes all tests in the BigCrush suite \cite{pierr}. For the largest of them the period is a number of nearly 
   a million digits. 
   The period q is defined in (\ref{period1}), (\ref{period}). }
   \label{tbl:largeN}
\end{table}

The spectrum of the operator $A(N,s)$   and of its inverse $A^{-1}(N,s)$
are presented in Figure \ref{fig2}.  The spectrum of the operator $A(N,s)$
has two real eigenvalues for even $N$ and three for odd $N$, all the rest
of the eigenvalues are complex and lying on leaf-shaped curves.
It is seen that the spectrum tends to a universal limiting form as
$N$ tends  to infinity, and the complex eigenvalues  $ 1/\lambda$
(of the inverse operator) lie asymptotically on the complex curve which has the representation
\be 
r(\phi) = 4 \cos^2(\phi/2)
\label{eq:curve}
\ee
in the polar coordinates $ \lambda = r \exp(i\phi)$. From the above 
analytical expression for eigenvalues it follows that the eigenvalues satisfying the 
condition $ 0 <  \vert \lambda_{\alpha} \vert   < 1$  are in the range $-2\pi/3 <  \phi < 2\pi/3$ 
and the ones satisfying the 
condition $ 1 <  \vert \lambda_{\beta} \vert $ are in the interval $2\pi/3 <  \phi < 4\pi/3$.
One can conjecture that there exists a limiting infinite-dimensional dynamical system with continuous space coordinate and discrete time with the above spectrum. The $A(N,s)$ system for a finite $N$ is then an approximation to this continuous dynamical system. For finite $N$, we found also the empirical formula for the eigenvalues of the $A(N,s=-1)$ system:
\be
\lambda_j = \frac{1}{4 \cos^2(j\pi/2N)} ~~ \exp(i\, \pi j/N) ~~ \textrm{for}~ j=-N/2..N/2 ~,
\label{eq:evs}
\ee
where the conjugate eigenvalues are for $\pm j$. This formula gives excellent approximation for the small eigenvalues of the operator $A(N,s)$ and is not applicable for the few of the largest eigenvalues. The derivation of the formulae \eqref{eq:curve}  and \eqref{eq:evs} will be given elsewhere.
The entropy of the C-system $A(N,s)$ can now be calculated for
large values of $N$  as an integral over eigenvalues (\ref{eq:curve}):
\be\label{linear}
h(A)= \sum_{\alpha   } \ln \vert {1 \over \lambda_{\alpha} }\vert = \sum_{-2\pi/3 <  \phi_i < 2\pi/3} \ln (4\cos^2(\phi_i/2)~ \rightarrow  ~N \int^{2\pi/3}_{-2\pi/3} \ln (4\cos^2(\phi/2){d\phi \over 2\pi}=
{2\over \pi}~ N
\ee
and to confirm  that it {\it increases linearly with the dimension $N$ of the operator} $A(N,s)$.

In the recent paper \cite{konstantin} the period of
the trajectories of the system $A(N,s)$ was found which is characterised by 
a prime number $p$\footnote{The
general theory of Galois field and the periods
of its elements can be found in \cite{mixmaxGalois,lnbook,nied,niki,alanen,pierre}.} .
In  \cite{konstantin} the necessary and sufficient criterion were
formulated for the sequence to be of the maximal possible period:
\be\label{period1}
q={p^N-1 \over p-1} \sim e^{(N-1)\ln p}.
\ee
{\it It follows then that the period of the trajectories exponentially increases with
the size of the  operator $A(N,s)$}. Thus the knowledge of the spectrum allows to
calculate the entropy (\ref{linear}) and the period (\ref{period1})   of the
trajectories.  The number of periodic trajectories (\ref{density}) behaves as
\be\label{density1}
\pi(q) \sim \exp{({2 N q \over \pi})} /q.
\ee
In summary, we found  the spectrum (\ref{eq:curve}), the entropy (\ref{linear}),
the period on a rational sublattice (\ref{period1}) and the corresponding  density (\ref{density})
of the C-system $A(N,s)$.
These quantities for large values of $N$ are presented in the Table \ref{tbl:largeN}
and for smaller values of $N$ the data can be found in \cite{konstantin}.

\section{\it Three-parameter Family of C-operators  $A(N,s,m)$ }

We note that the special form of the matrix $A(N,s)$ in \eqref{eq:matrix} has the highly desirable property of having a widely spread, nearly continuum spectrum of eigenvalues (see Fig.\ref{fig2}),
which indicates that the mixing of the dynamical system is occurring on all scales \cite{yer1986a}.
This property appears to be a consequence of its very special, near-band-matrix form. At the same time, the last column assures that the determinant of the matrix is equal to one, and therefore the phase volume of the dynamical system is conserved and fulfils 
the conditions (\ref{mmatrix}). 

The construction of the matrix $A(N,s)$ uses 
increasing natural numbers.  On the last line of the matrix $A(N,s)$ in 
\eqref{eq:matrix} we used the natural numbers 
from 1 to $N$. This means that the entries of the matrix $A(N,s)$ are gradually 
increasing with $N$  and it results in the corresponding increase of the eigenvalues 
\eqref{eq:curve}  of the operator $A(N,s)$ and of its entropy $h(A)$, as in \eqref{linear} . 
The other advantage of the special form of the matrix $A(N,s)$ is that it allows 
to generate a computer code which executes fast multiplication of the matrix $A(N,s)$ 
and with the vector. Our aim is to find out the generalisations of the operator $A(N,s)$ which 
keep intact all the above properties and generate even larger entropies at a given 
size  $N$ of the operators.

A three-parameter family of  operators $A(N,s,m)$, which we present next, is constructed by
replacing the sequence in the bands, below the diagonal, which is originally $3,4,5,...,N$ with the sequence $m+2,2m+2, 3m+2, ...,(N-2)m+2$, where $m$ is some integer:
\be
\label{eq:matrix1}
A(N,s,m) =
   \begin{pmatrix}
      1 & 1 & 1 & 1 & ... &1& 1 \\
      1 & 2 & 1 & 1 & ... &1& 1 \\
      1 & m+2+s & 2 & 1 & ... &1& 1 \\
      1 & 2m+2 & m+2 & 2 &   ... &1& 1 \\
      1 & 3m+2 & 2m+2 & m+2 &   ... &1& 1 \\
      &&&...&&&\\
      1 & (N-2)m+2 & (N-3)m+2 &  (N-4)m+2 & ... & m+2 & 2
   \end{pmatrix}
\ee
\begin{figure}
\begin{center}
\includegraphics[width=3cm]{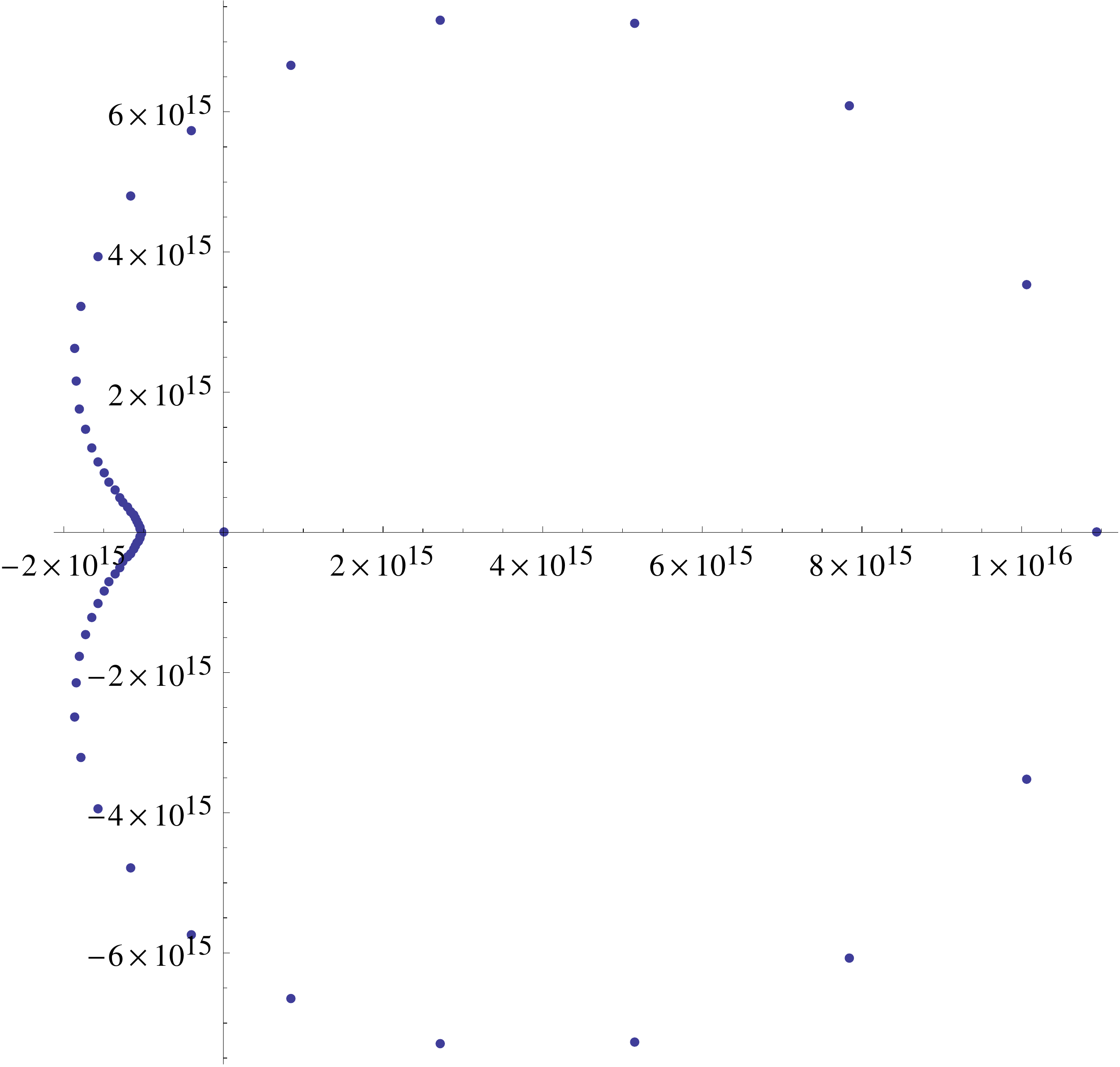}~~~~~~~
\includegraphics[width=3cm]{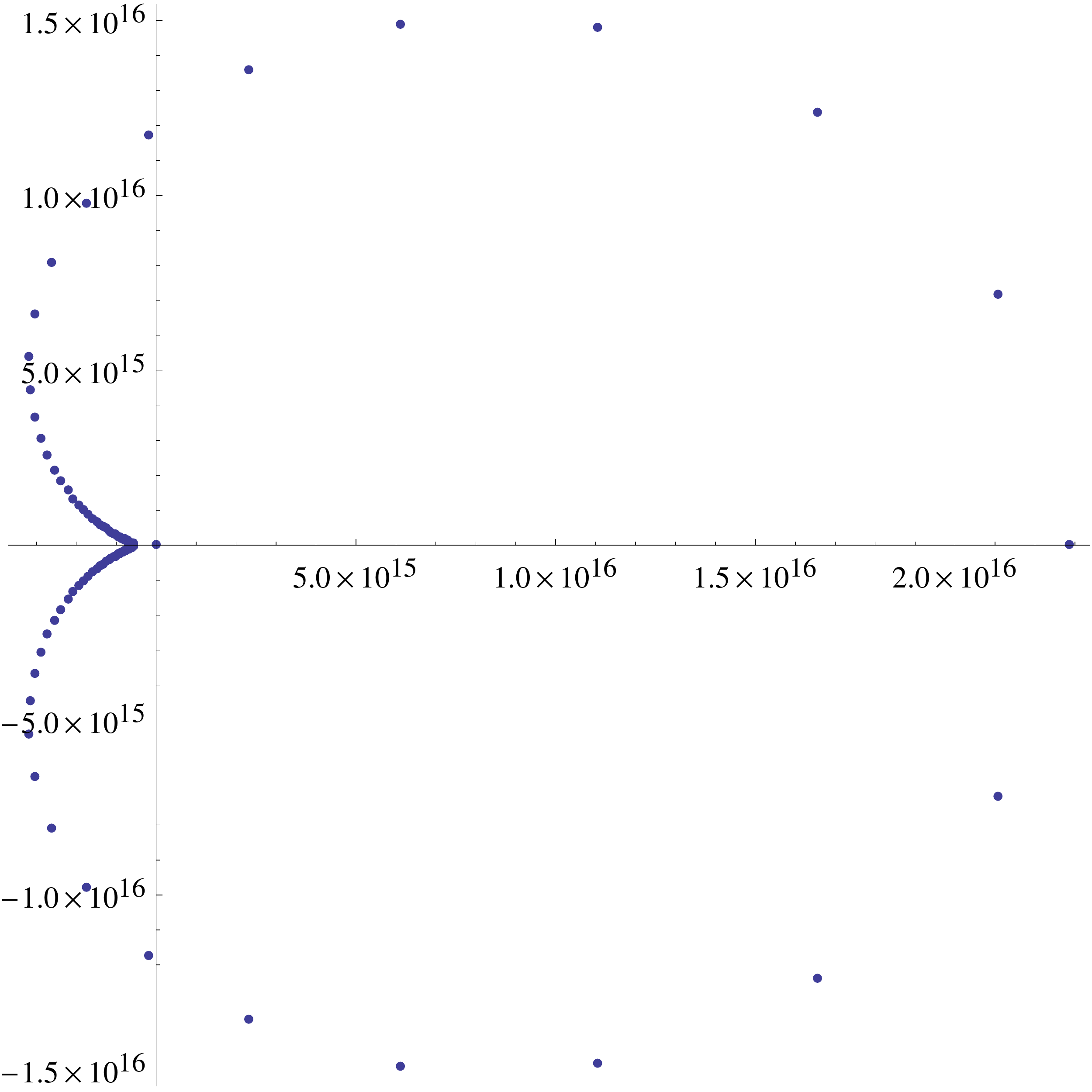}~~~~~~~
\includegraphics[width=3cm]{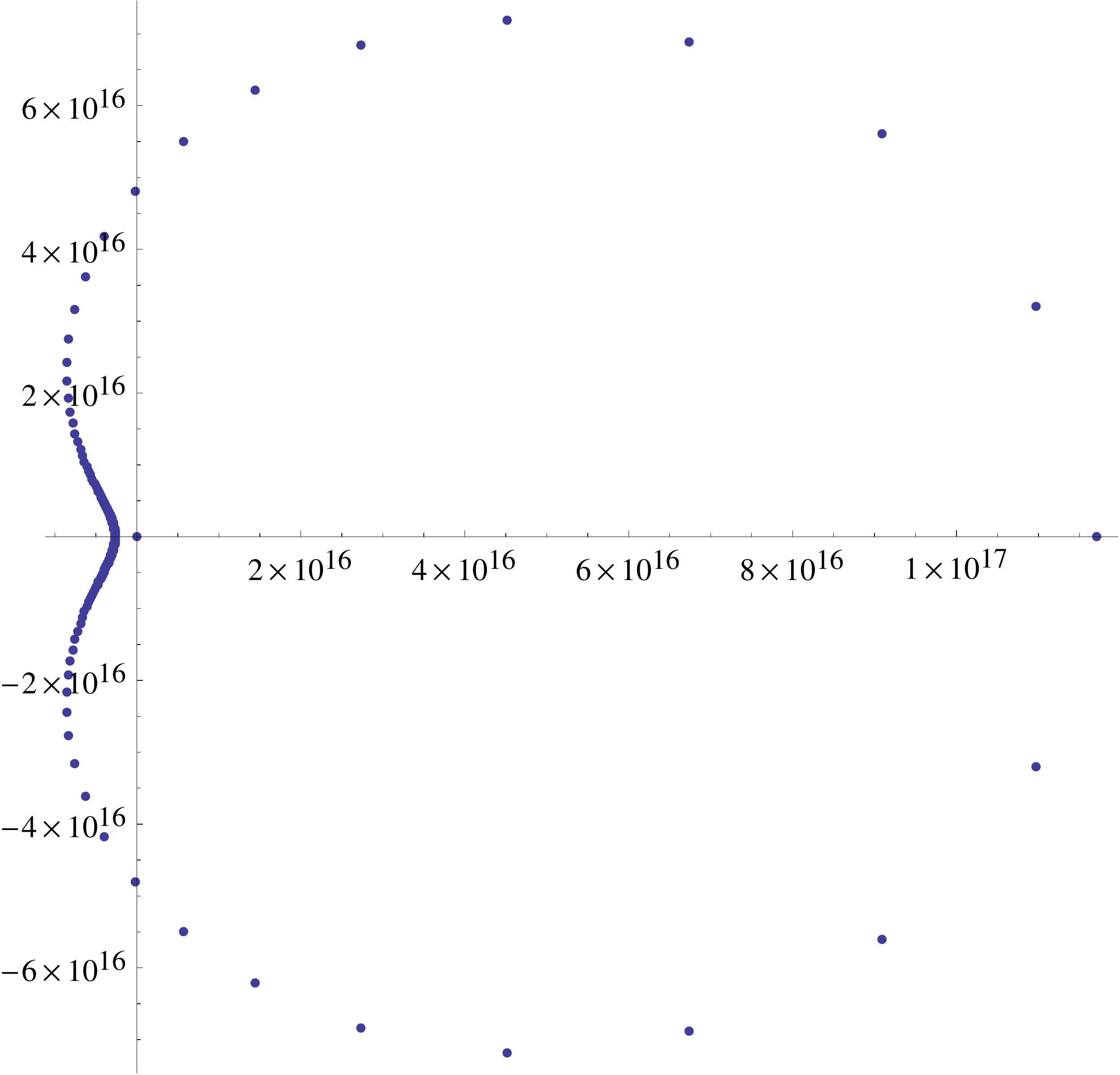}
\caption{The distribution of the eigenvalues of the operator $A(N,s,m)$  in (\ref{eq:matrix1})
for large value of $s$ and $m$, but fixed, and increasing values of $N=60, 120, 240$
from Table \ref{tbl:largeM}.
The spectrum represents  a leaf  of a large radius
proportional to $ \lambda_{max} \approx m$
and a very small eigenvalue at the origin  $ \lambda_{min}\approx m^{-N+1}$.
With  increasing $N$ the ``stem of the leaf" becomes more pronounced on the left hand side
of the spectral curve. }
\label{fig3}
\end{center}
\end{figure}
Thus the case of $m=1$ simply corresponds to the original matrix (\ref{eq:matrix}). It is most advantageous to take large values of $m$, but preferably keeping $Nm < p$, such as to have an unambiguous correspondence between the continuous system \eqref{eq:rec} and the discrete system on the rational sublattice.  The distribution of the eigenvalues of the operator $A(N,s,m)$
for the values of s and m which are given in Table \ref{tbl:largeM} are presented
on Fig.\ref{fig3}. The spectrum for the increasing values of $N=60, 120, 240$ is shown
in sequence, from left to right.
The spectrum represent a leaf  of a large radius proportional
to $\lambda_{max} \approx m$
and very small eigenvalue at the origin  $ \lambda_{min}\approx m^{-N+1}$.
With  increasing $N$ the ``stem of the leaf" becomes more pronounced on the left hand side
of the spectral curve.

\begin{table}[htbp]
   \centering
   \begin{tabular}{@{} lcccccl @{}}
      \toprule
      Size & Magic &Magic&Entropy& Log of the period  q\\
      N    & $m$ & s& &$\approx \log_{10} (q)$  \\
      \midrule
     8    & $m=2^{53}+1$&s=0 &220.4& 129\\
      17  & $m=2^{36}+1$&s=0& 374.3&294\\
      40  & $m=2^{42}+1$&s=0&1106.3& 716\\
      60  & $m=2^{52}+1$&s=0&2090.5& 1083\\
      96  & $m=2^{55}+1$&s=0&3583.6&1745 \\
      120  & $m=2^{51}+1$ & s=1&4171.4& 2185\\
       240  & $m=2^{51}+1$& s=487013230256099140&8679.2& 4389\\

    \bottomrule
   \end{tabular}
   \caption{Table of three-parameter MIXMAX generators $A(N,s,m)$ in (\ref{eq:matrix1}). These generators have an advantage of having a very high quality sequence for moderate and small $N$. In particular, the smallest generator we tested, $N=8$, passes all tests in the BigCrush suite
   \cite{pierr}.   The period $q$ is defined in (\ref{period1}), (\ref{period}). }
   \label{tbl:largeM}
\end{table}

A further possible generalization of the three-parameter family of operators is the following:
the four-parameter  operators $A(N,s,m,b)$ is constructed by
replacing the sequence in the bands, below the diagonal, which is originally $3,4,5,...,N$ with the sequence $3m+b, 4m+b, 5m+b,...,Nm+b$, where $m$ and $b$ are some integers:
\be
\label{eq:matrix2}
A(N,s,m,b) =
   \begin{pmatrix}
      1 & 1 & 1 & 1 & ... &1& 1 \\
      1 & 2 & 1 & 1 & ... &1& 1 \\
      1 & 3m+s+b & 2 & 1 & ... &1& 1 \\
      1 & 4m+b & 3m+b & 2 &   ... &1& 1 \\
      1 & 5m+b & 4m+b & 3m+b &   ... &1& 1 \\
      &&&...&&&\\
      1 & Nm+b & (N-1)m+b &  (N-2)m+b & ... & 3m+b & 2
   \end{pmatrix} .
\ee
This four-parameter family $A(N,s,m,b)$ reduces back to the three-parameter family $A(N,s,m)$ for $b=2-2m$.
It is the case that some of these four-parameter generators, for specially chosen $m$ and $b$,  allow efficient computer multiplication -
the property which plays an essential role if one tries to use these operators for Monte-Carlo
simulations.

\section{\it Computer Implementation. MIXMAX Random Number Generator}

In a typical computer implementation \cite{konstantin,hepforge} of the periodic trajectories 
of the automorphism \eqref{eq:rec} the initial vector 
$$
u^{(k)}= (u^{(k)}_1,...,u^{(k)}_N)
$$ has a rational
components $u^{(k)}_i=a^{(k)}_i/p$, where $a^{(k)}_i$ and $p$ are natural numbers
and $i=1,...,N$.
Therefore it is convenient to represent $u^{(k)}_i$ by its numerator $a^{(k)}_i$ in
computer memory and define the iteration in terms of $a^{(k)}= (a^{(k)}_1,...,a^{(k)}_N)$
\cite{konstantin,mixmaxGalois}:
\be
\label{eq:recP}
a^{(k+1)}_i = \sum_{j=1}^N A_{ij} \, a^{(k)}_j ~~~\textrm{mod}~ p~, ~~~~~~~k=0,1,2,....
\ee
If the denominator $p$ is taken to be a prime number \cite{konstantin,mixmaxGalois,hepforge},
then the recursion is realised on extended
Galois field $GF[p^N]$  \cite{nied,niki} and it
allows to find the period of the trajectory $q$ in terms of $p$ and the properties of the
characteristic polynomial $P(x)$ of the matrix $A$ \cite{konstantin,mixmaxGalois,alanen,hepforge}. If
the characteristic polynomial $P(x)$ of some matrix $A$ is primitive in the
extended Galois  field $GF[p^N]$, then
\cite{mixmaxGalois,lnbook,nied,alanen,pierre}:
\be\label{period}
 A^q = p_0~ \mathbb{I}~~\textrm{ where}~~  q=\frac{p^N-1} {p-1} ~, 
 \ee
where $p_0$ is a free term of the  polynomial $ P(x)$ and is a {\it primitive element} of $GF[p]$.
Since our matrix $A$ has $p_0=\Det  A= 1$, the polynomial $ P(x)$ of $A$ cannot be primitive.
The solution suggested  in \cite{konstantin} is to define the necessary and
sufficient conditions for the period $q$
to attain its maximum, and they are \cite{konstantin}:
\begin{enumerate}
\item[\bf{1.}] $A^q = \mathbb{I} ~(mod~ p) $,~~~where $q=\frac{p^N-1} {p-1}$.
\item[\bf{2.}] $A^{q/r} \neq \mathbb{I} ~(mod~ p)$,~~~~ for any $r$ which is a prime divisor of $q$ .
\end{enumerate}
The first condition is equivalent to the requirement  that the characteristic polynomial is {\it irreducible}.
The second condition can be checked if the integer factorisation of $q$ is available \cite{konstantin}, then
the period of the sequence is equal to (\ref{period}) and is independent of the seed.
There are precisely $p-1$ distinct  trajectories which together fill up all states of the $GF[p^N]$ lattice:
\be
 q~ (p-1) = p^N-1.
\ee
In \cite{konstantin} the actual value of $p$ was  taken as $p=2^{61}-1$,
the largest Mersenne number that fits into an
unsigned integer on current 64-bit computer architectures.
For the purposes of generating pseudo-random numbers
with this method, one chooses any initial vector $a^{(0)}$,  called the "seed",
with at least one non-zero component.

\begin{figure}
\begin{center}
\includegraphics[width=5cm]{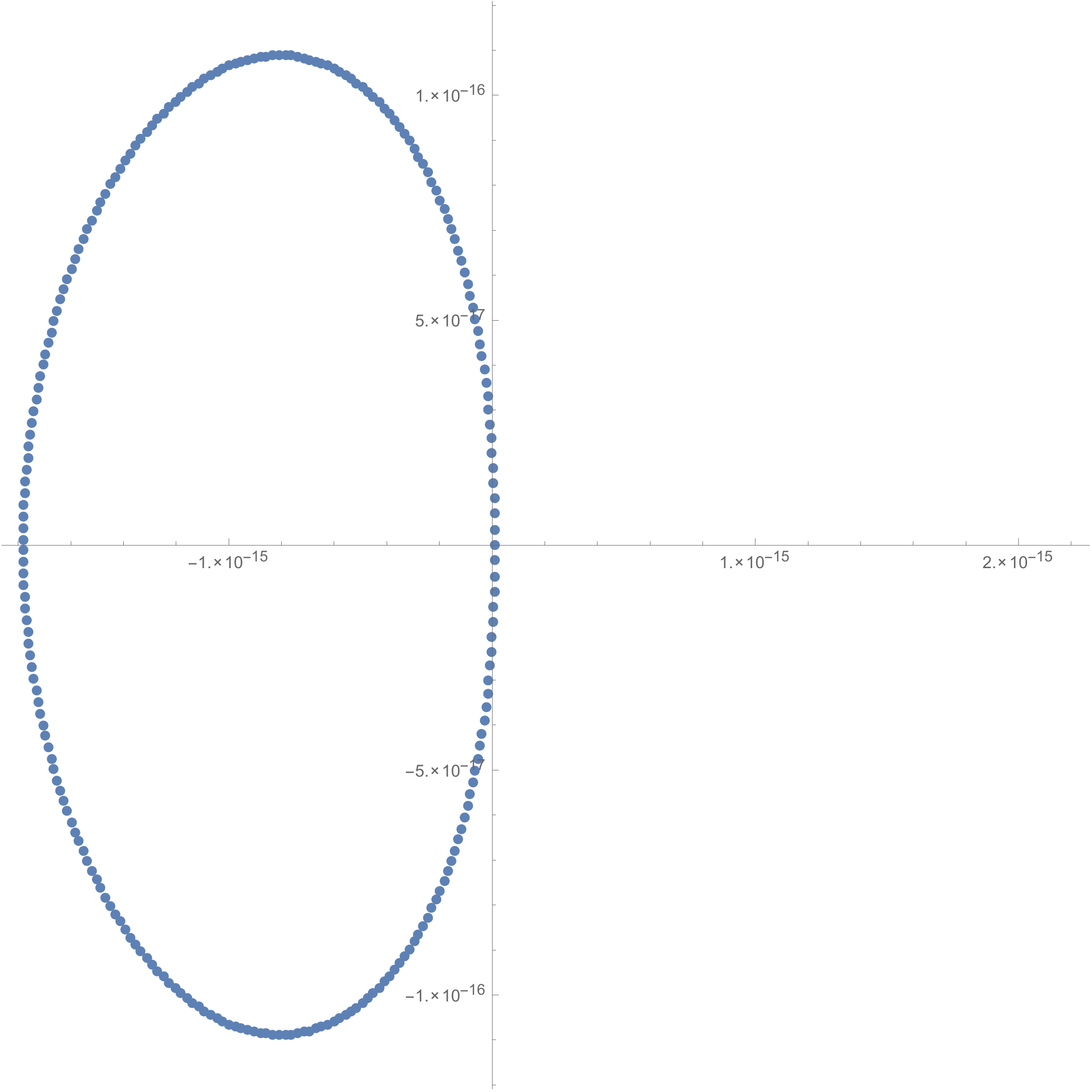}
\caption{The distribution of the eigenvalues of the inverse operator $A^{-1}(N,s,m)$
of (\ref{eq:matrix1}) for large values $N=240,  s=487013230256099140$ and $m=2^{51}+1$. }
\label{fig4}
\end{center}
\end{figure}

The commonly-used RNGs based on a 
linear recurrence typically reach the maximal period of $ p^N-1$ using the primitive 
elements of $GF[p^N]$. 
They use either a large prime number $p$ as in \cite{pierr,pierre,grothe,knuth1,pierre3}, 
or $p = 2$ as in \cite{pierre5,matsumoto,nishimoura,panted}.

The results of the last two sections allow to disclose some additional parameter values
for the MIXMAX generator,
in addition to those found in \cite{konstantin}.
The properties of the MIXMAX generators improve appreciably with $N$, the size of the matrix,
and therefore we undertook a search for large values of $N$ and some small values
of the parameter $s$.
Because the speed of the generator does not depend on $N$,
these generators are useful if the dimension $D$ of the Monte-Carlo integration is large but finite,
in which case one would like to choose $N \geq D$. If a generator with such large $N$ is available,
then the convergence of the Monte-Carlo result to the correct value and with a residual
which is normally distributed is assured. The latter guarantee is given by the theorem of Leonov
\cite{leonov,Savvidy:2015ida}.

Our search of the parameters for the MIXMAX generator $A(N,s)$
with large $N$ and maximal period  has yielded the values presented in the Table \ref{tbl:largeN}.
As one can deduce from this data, the entropy is
linearly increasing with $N$ (\ref{linear}). As it was demonstrated in \cite{konstantin}, the Kolmogorov entropy, which needs to be greater than about $h(A) \approx 50$ for the generator to be empirically acceptable.
Therefore it should not be surprising that for all of these generators, the sequence  passes all tests in the BigCrush suite \cite{pierr}. For the largest of them $N=44851$, the period is a 
number of nearly a million digits.  If an increase in entropy is desired without increasing the size of the matrix $N$, it is
now possible to search for large $s$ as well.
The combinations of $N$ and $s$ which we found to be useful in this regard are
presented in Table \ref{tbl:largeS}.
The generator with $N=256$ and $s=487013230256099064$ has the best combination of speed, reasonable size of the state, and availability for implementing the parallelization by skipping and is currently available  generator in the ROOT and CLHEP software packages at CERN
 for scientific calculation \cite{cern,root,clhep}.

For the  three-parameter family $A(N,s,m)$ of the MIXMAX generators the
 convenient values of the parameters are provided in Table \ref{tbl:largeM}.
The efficient implementation in software can be achieved for some particularly
convenient values of $m$ of the form $m=2^k+1$ \cite{hepforge}.
Inspecting the data in the Table \ref{tbl:largeM} one can get convinced that
the system with $N=240,  s=487013230256099140$ and $m=2^{51}+1$
has the best stochastic properties within the $A(N,s,m)$ family of
operators.  It is also obvious that the record value of the
entropy $h(A)= 28704.6$ which we achieved for the  generator $A(N,s)$
with $N=44851,~s=-3$ remains the biggest and required few months of 
computer search.

\section{\it Acknowledgement }
This work was supported in part by the European Union's Horizon 2020
research and innovation programme under the Marie Sk\'lodowska-Curie
Grant Agreement No 644121.

\section{\it Note Added }
In this note we shall provide some additional comments. 

The phase space of the  C-systems under consideration (\ref{eq:rec} )
is a $N$-dimensional torus \cite{anosov,yer1986a,konstantin,Savvidy:2015ida}, 
appearing at  factorisation of the 
Euclidean space $E^N$ with coordinates $u= (u_1,...,u_N)$ over an integer lattice $\CZ^N$. 
The coordinates 
$u$ on a torus are rational and  irrational numbers. The operator $A$ acts on the 
initial vector $u^{(0)}$ and produces a phase space trajectory $u^{(n)}=  A^n u^{(0)}$
on a torus.  The trajectories 
of a C-system can be periodic and non-periodic. {\it  All trajectories 
which start from vectors $u^{(0)}$ which have rational coordinates, and only 
they, are periodic} \cite{anosov,Savvidy:2015ida}.  
The rational numbers are everywhere dense on the phase space of a 
torus and {\it the periodic trajectories of the C-systems follow the same pattern} 
and are everywhere dense  \cite{anosov}, 
like rational numbers on a real line. 

It is important now to raise the following fundamental question: 
{\it  can one approximate, and approximate uniformly,  
non-periodic trajectories of a C-system by  using periodic trajectories of the same 
C-system}?    The answer is yes \cite{anosov}.

The trajectories which start from the irrational 
coordinates are uniformly distributed over 
the phase space on a torus and the main goal is to approximate, and approximate 
uniformly, these non-periodic trajectories.
Therefore what we suggested in our preprint articles which were released in 1986 
\cite{yer1986a} and finally published  in 1991 is
to {\it use  the periodic trajectories of a  C-system 
to approximate, as well as possible, the non-periodic, uniformly distributed 
trajectories of the same C-system}. The periodic and non-periodic trajectories  of 
a C-system should not be separated from each other 
within the C-system, they are in a harmonic union, as the rational and irrational 
numbers on a real line are. 

We shall provide an example from calculus which may help 
to understand this harmonic union of trajectories.  
A good illustration example will be the approximation of the non-algebraic, transcendental 
number $\pi$ by rational numbers, which was the aim of  Archimedes, 
$3{1\over 7} < \pi < 3{10\over 71}$, of Ramanujan  and the series of Leibniz: 
$$
\pi = 4 \sum^{\infty}_{n=0}{(-1)^n\over 2n+1} =4(1 - {1\over 3}+ {1\over 5}-..... )~.
$$
One should ask: is  this representation of $\pi$ in terms of rational numbers useful to us, 
humans, who are unable to 
embrace the number $\pi$ on a computer? 
It seems that it is  a unique and  useful   way  to 
grasp  the number $\pi$, and, indeed, 13.3 trillion digits of $\pi$ were computed in October 2014. 
The same is true for the C-systems, the aim is to 
"follow" the non-periodic trajectories {\it for as long and 
as closely}  as possible by the periodic trajectories. The analogy will be  complete if one 
recalls  that the $\pi=3.14~....$ expressed in terms of decimal digits, is a non-periodic 
number, while the rational numbers are periodic. In this case we  have  approximation of 
a non-periodic number by a periodic. 

It is therefore important to learn  how dense are the periodic 
trajectories populating the phase space of a C-system and how this density is 
correlated with the entropy.
The formula (\ref{density})  describes the asymptotic 
 distribution of the periodic trajectories of the length $q$,  quantifying the rate at which they occur (see \cite{Savvidy:2015ida} and references therein).  The asymptotic formula (\ref{density}) approximates $\pi(q)$ in the sense that the relative error of this approximation approaches zero as $q$  tends to infinity, but does not say anything about the limit of the difference $$\pi(q) - {e^{q \ h(A)} \over  q}. $$   
 A good analogy will be the prime number theorem. There are infinitely many primes, as demonstrated by Euclid around 300 BC and the prime number theorem 
 describes the asymptotic distribution of the prime numbers among the positive integers.

 In this article we calculated some of 
 the periods and it is  obviously 
a nontrivial task.  The Table 2 which presents the specific values of $N$ and $s$ and
requires development of a computer code and months of computer time, the task 
which is analogous 
of searching large prime numbers.
For the  MIXMAX generator one should find those values 
of $N$ and $s$ for which all the necessary conditions are fulfilled. 
The fact that here we reached a period which is a million digits 
was not a "sport" achievement but the ability to "follow" a non-periodic trajectory as 
long as possible. The ability to increase further the dimension $N$, the entropy 
 and the periods is a priceless advantage of the MIXMAX family  of RNGs allowing 
 them, in principle,  and it is justified theoretically, to pass even stronger tests
 \cite{marwan,mice}.

\vfill

\end{document}